\documentclass[12pt]{article}
\usepackage{amsmath, amssymb, latexsym, amscd, amsthm, enumitem}

\def\Int{\operatorname{Int}}
\def\id{\operatorname{id}}
\def\diam{\operatorname{diam}}

\def\Int{\operatorname{Int}}
\def\id{\operatorname{id}}
\def\dem{\operatorname{dem}}
\def\diam{\operatorname{diam}}
\def\Id{\operatorname{Id}}
\def\Homeo{\operatorname{Homeo}}
\def\Isot{\operatorname{Isot}}
\def\Int{\operatorname{Int}}

\newtheorem{theorem}{Theorem}[section]
\newtheorem{corollary}[theorem]{Corollary}

\newtheorem{proposition}[theorem]{Proposition}

\theoremstyle{definition}
\newtheorem{definition}[theorem]{Definition}

\begin{document}

\title{ON PROJECTIONS OF A COMPACT SET IN $\mathbb R^N$\\
UNDER THE ACTION OF A TYPICAL AMBIENT HOMEOMORPHISM} 
\author{OLGA FROLKINA\\
Chair of General Topology and Geometry\\
Faculty of Mechanics and Mathematics\\
M.V.~Lomonosov Moscow State University\\
Leninskie Gory 1, GSP-1,\\
Moscow 119991, Russia, \\
olga-frolkina@yandex.ru
}

\date{}

\maketitle

\renewcommand{\thefootnote}{}



\centerline{received (30 May 2026)}


\bigskip

\maketitle

\begin{abstract}
We apply
ideas of geometric measure theory and
Baire category theory 
to topological problems,
namely, to topological embeddings
of compact sets into Euclidean spaces.

In 1947, Borsuk constructed a Cantor set
in $\mathbb R^N$, $N\geqslant 3$, such
that its projection onto any $(N-1)$-plane
contains an $(N-1)$-dimensional ball.
This can be strengthened:
a desired Cantor set 
can be obtained
from an arbitrary Cantor set by an
arbitrarily small isotopy of the space $\mathbb R^N$.
The question arises: how do the dimensions of the projections of a compact set
$X\subset \mathbb R^N$ behave under a typical ambient 
isotopy or under a typical ambient homeomorphism?
(Typical in the sense of the Baire category.)
We solve this problem.
As a consequence, we get new criteria
of tameness and wildness of a Cantor set
in terms of its projections.
Our main result
strengthens
V{\"a}is\"{a}l\"{a}'s theorem (1979) connecting 
Hausdorff dimension
 and Shtan'ko embedding dimension.
In its turn, V{\"a}is\"{a}l\"{a}'s theorem
extends results of
N\"{o}beling (1931) and Szpilrajn (1937)
on relationship between
Hausdorff dimension and
topological dimension.
\end{abstract}

\section{Introduction}\label{introduction}

L.~Antoine constructed a Cantor
set in plane such that all of its projections coincide with those of a regular hexagon
\cite[{\bf 9}, p. 272--273]{Antoine-FM}.
Other examples with similar properties were
described by Otto, A.~Flores,
G.~N\"{o}beling
\cite{Ergebnisse}; see also \cite{Dijkstra-van-Mill},
\cite[Stat. 11]{Frolkina-Greece-volume}.

For each $N\geqslant 2$,
K.~Borsuk
constructed a Cantor set
in $\mathbb R^N$ such that its projection into any $(N-1)$-plane
contains a ball,
equivalently, has dimension $(N-1)$ \cite{Borsuk}.
See also \cite[Thm. 6.2]{Bing-tame}.
Many papers strengthen and extend
Borsuk's result, see e.g.
\cite{Meyerson},
\cite{Glaser},
\cite{Kuzminykh},
\cite{Halverson-Wright},
\cite{Goodsell-Wright}.

In \cite{Cobb-projections}, J.~Cobb described
a Cantor set in $\mathbb R^3$
such that its projection into any $2$-plane
is $1$-dimensional, and posed a general question:                  
for given integers
$N>m\geqslant k>0$, 
does there exist
a Cantor set in $\mathbb R^N$ such that its projection
into any 
$m$-plane has topological dimension~$k$~?
(We call these sets briefly $(N,m,k)$-sets.)
The case of $m=k$ is answered positively by Borsuk's construction.
For the cases
$(N,m,k=m-1)$ and $(N, m=N-1 , k)$, desired sets 
were described in
\cite{Frolkina-proj} and \cite{BDM}, correspondingly; 
both of these papers extend Cobb's method.
Applying facts from the theory of tame and wild
Cantor sets, Frolkina obtained new
wide series of  $(N,N-1,N-1)$- and $(N,N-1,N-2)$-sets 
\cite{Frolkina-Greece-volume},
\cite{Frolkina2021}.

\begin{definition}
A zero-dimensional compactum $X\subset \mathbb R^N$ 
is called \emph{tame}
if 
there is a homeomorphism
$h:\mathbb R^N\cong \mathbb R^N$
such that
$h(X)$ is a subset of a straight line.
Otherwise, $X$ is called \emph{wild}.
\end{definition}

Any zero-dimensional compactum in plane is tame
\cite[{\bf 75}, p.~87--89]{Antoine-diss},
see also
\cite[Cor. II.3.2, II.3.3]{Keldysh} or
\cite[Chap.~13]{Moise}.

The first wild Cantor sets
in $\mathbb R^3$ 
were described by L.~Antoine
in 1920--21
and by P.S.~Urysohn in 1922--1923 (independently).
Antoine's construction was later extended to
 $\mathbb R^N$, $N\geqslant 3$;
now there are essentially different examples with additional properties.
For references, see \cite{Frolkina-fences}.

Using his $(N,N-1,N-1)$-set,
Borsuk constructed a simple arc in $\mathbb R^N$ such that
its projection into any $(N-1)$-plane contains an
$(N-1)$-ball \cite[p. 277]{Borsuk}. 
The same argument provides a knot in $\mathbb R^N$ with the same property.
For $N=3$ this answers a question of R.~Fox.
Moreover: analyzing Borsuk's paper, 
we see that a knot with this property exists in any equivalence class of knots 
(no matter, tame or wild).
See Proposition \ref{intro}.

\begin{definition}\label{tame-P}
A subset of $\mathbb R^N$
is called a \emph{polyhedron}
if it is the union of a finite collection of simplices.
A compactum $X \subset\mathbb  R^N$  
homeomorphic to a polyhedron is
called
\emph{tame}
if there exists a homeomorphism  
$h$ of $\mathbb R^N$ onto itself 
such that $h(X)$ is a polyhedron in $\mathbb R^N$;
otherwise, $X$ is \emph{wild}.
\end{definition}

Compacta which are embedded in a 
topologically equivalent way
may behave absolutely dif\-fer\-ent
if we consider dimensions of their projections.
Proposition~\ref{intro} includes
several statements of this type.
(For a tame Cantor set, compare 1a) and 1b).
For a wild Cantor set, compare 2) and 3).
For a knot, see 3) and 4).)

\begin{definition}\label{isotopy-def}
An \emph{isotopy} of a space $Y$
is%
\footnote{Here and below, $I=[0,1]$.}
 a continuous map $F  :Y\times I \to Y$ such that
$f_t $ is a homeomorphism $Y\cong Y$ for each $t\in I$,
and 
$f_0 = \id $.
(As usual, $f_t = F|_{Y\times \{t\}}  $. 
An isotopy $F$ will be also denoted as $\{f_t \} : Y\cong Y$.)

An isotopy $\{ f_t\} : Y\cong Y$
is called an
\emph{$\varepsilon $-isotopy}
if
$d(x,f_t(x)) \leqslant \varepsilon $
for each $t\in I$ and each
$x\in Y$.

The support of an isotopy
 $\{ f_t\} : Y\cong Y$ is the closure of the set
$$
\{ x\in Y \ | \ \text{for some } t\in I \text{ we have } x\neq f_t(x)  \} .
$$
\end{definition}

\begin{proposition}\label{intro}
Let $N\geqslant 2$ be an integer, and $\varepsilon >0$.

1) Let $X\subset \mathbb R^N$ be a tame Cantor set. 
Then%
\footnote{Here and below, $P_\Pi (M)$
denotes the orthogonal projection of a set
$M\subset \mathbb R^N$ into a linear subspace
$\Pi\subset \mathbb R^N$.
As usual, $O_{\varepsilon }X $ is the open $\varepsilon $-neighborhood of $X$.
}
\\
a) 
there exists an $\varepsilon $-isotoply $\{ h_t \} : \mathbb R^N\cong \mathbb R^N$ with support in $O_\varepsilon X$
such that
$\dim P_{\Pi } (h_1 (A)) = \dim \Pi $ 
for any non-empty open subset $A\subset X$ and
any proper linear subspace
$\Pi \subset \mathbb R^N$;
\\
b) 
there exists an
$\varepsilon $-isotopy $\{ h_t \} : \mathbb R^N \cong \mathbb R^N$
with support in $O_\varepsilon X$
such that
the set $h_1(X)$ has general position with respect to all projections 
\cite[Def. 3.1]{Frolkina-arxiv}.
In particular, $P_\Pi (h_1(X))$ is a Cantor set
for any non-zero linear subspace
$\Pi \subset \mathbb R^N$.

2) Let $X\subset \mathbb R^N$ be any Cantor set.
There exists 
an
$\varepsilon $-isotopy $\{ h_t \} : \mathbb R^N\cong \mathbb R^N$
with support in
 $O_\varepsilon X$
 such that
$h_1(X)$
is an ${(N,N-1,N-2)}$-set.

3)
Let $X\subset \mathbb R^N$ be an uncountable compactum.
There exists 
an 
$\varepsilon $-isotopy $\{ h_t \} : \mathbb R^N\cong \mathbb R^N$
with support in $O_\varepsilon X$
such that
$\dim P_{\Pi } (h_1 (X)) = \dim \Pi $
for any proper linear subspace
$\Pi \subset \mathbb R^N$.

4) Suppose that $X\subset \mathbb R^3$ 
is homeomorphic to a circle or to a segment.
There exists an
$\varepsilon $-isotopy $\{ h_t \} : \mathbb R^3 \cong \mathbb R^3$
with support in $O_\varepsilon X$
such that
${\dim P_\Pi (h_1(X)) = 1}$
for any non-zero linear subspace
$\Pi \subset \mathbb R^3$.
\end{proposition}

The aim of this paper is
\emph{
to understand
which behavior of isotopies
is typical for a given compactum $X\subset\mathbb R^N$.}
``Typical'' is understood in the sense of Baire category.
Below we answer this question.
At the same time, we strengthen  
V{\"a}is\"{a}l\"{a}'s theorem \cite{Vaisala}.

Formally, Proposition~\ref{intro} is new.
Its parts 1a), 3) strengthen the result from \cite{Borsuk}
and from 
\cite[Stat. 12]{Frolkina-Greece-volume}.
Statement 1b) implies
\cite[Thm. 1]{Frolkina-Arch-volume},
\cite[Thm. 5]{Cobb-projections}.
Part 2) implies \cite[Thm. 2]{Frolkina-Greece-volume}.
Part 4) implies \cite{Borsuk}.
But main ideas which prove Proposition~\ref{intro}
are contained in
\cite{Borsuk}, 
\cite{Frolkina-Greece-volume}.
Proposition~\ref{intro}
is stated here to motivate the main problem.
We prove it in Section~\ref{intro-proof}.
Part 4) is strengthened in Corollary~\ref{case-knots}.

Similar theorems about projections exist in measure theory:
instead of to\-po\-lo\-gi\-cal dimension, the Hausdorff dimension is considered, and projections are taken
on ``almost all" planes,
i.e., with the exception of the set of planes with measure zero.
See \cite[Cor.~9.4, Thm. 16.2, 18.1]{Mattila},
\cite[Ch.~6]{Falconer}, \cite{Dijkstra-van-Mill}.

\subsection{Agreements, notation}

$\overline A$ and
$\partial A$ are the closure and the boundary of a set $A$, correspondingly.

$\dim X$ denotes topological dimension of a separable
metric space $X$ \cite{HW}, \cite{Engelking}
(for such spaces, three classical topological definitions of dimension
are equivalent
\cite[Thm. 4.1.5]{Engelking}).

$m_q (X)$ is the $q$-dimensional Hausdorff measure,
and
$\dim _H X$ is the Hausdorff dimension of a space $X$. See \cite{HW}.
 
 A linear subspace
 $\Pi \subset \mathbb R^N$
 is \emph{proper} if
 $\{ 0 \} \neq \Pi  \neq \mathbb R^N$.

\section{About the concepts needed to formulate an answer}

\subsection{Embedding dimension}

We need Shtan'ko's ``embedding dimension" theory.

\begin{definition}\label{isotopy-def}
A \emph{pseudoisotopy} of a space $Y$
is a continuous map $F  :Y\times I \to Y$ such that
$f_t $ is a homeomorphism $Y\cong Y$ for each $t\in [0,1)$,
and 
$f_0 = \id $.
(That is, $h_1$ need not be a homeomorphism.)
\end{definition}

\begin{definition}\cite{Stanko1970}
For a non-empty compactum
$X\subset\mathbb R^N$,
\emph{the dimension of embedding}
is the smallest integer
$\dem X = k \geqslant 0$
with the property:
for each $\varepsilon >0$
there is an $\varepsilon $-pseudoisotopy
$\{ f_t  \} : \mathbb R^N\to \mathbb R^N$
such that
its support lies in $O_\varepsilon X$
and
$f_1 (X)$ is a $k$-dimensional polyhedron in 
$\mathbb R^N$.
\end{definition}

The notation ``dem'' is an abbreviation
for ``dimension of embedding''.
We have $\dim X \leqslant \dem X$.
This number is not a topological invariant of $X$,
however,
it is an invariant of the equivalence class of the embedding
$X\hookrightarrow \mathbb R^N$. This means that
$\dem X = \dem h(X)$ for any homeomorphism
$h:\mathbb R^N\cong\mathbb R^N$ \cite[Thm.~3]{Stanko1970}.
Briefly, for $\dem X = k$, the embedded compactum
$X$ behaves ``like'' a $k$-dimensional polyhedron.
This concept was introduced by Stan'ko with the aim of extending
the definition of tame and wild embeddings to arbitrary compact sets. Previously, definitions of tameness and wildness were available only for
embeddings of
polyhedra and zero-dimensional compacta.

For a Cantor set $K\subset \mathbb R^N$, 
tameness is equivalent to the equality
$\dem K = 0$
\cite[Remark 3]{Stanko1970}, \cite[Thm.~I.4.2]{Keldysh}.

The 
demension of each wild Cantor set
in $\mathbb R^N$
equals $N-2$
\cite{Bryant1969}, 
\cite{Stanko1970},
\cite[Thm. 1.4]{Edwards}, \cite[Thm. 3.4.11]{DV}
(for $N=4$, see also
\cite{Shtanko-4dim}, \cite[p.~5]{BDVW}, \cite[Thm.~2.5.1]{Quinn},
\cite{Freedman}, \cite{Freedman-Quinn}).

The case of $1$-dimensional compacta in $\mathbb R^3$ is special. 
In particular, any embedding of a segment or circle in $\mathbb R^3$ is tame in the sense of Shtan'ko, i.e., its embedding dimension 
equals~$1$~\cite[Satz~4]{Bothe-U} .

For details
(in particular, other equivalent
definitions of the dimension $\dem $) see \cite{Stanko1970}, \cite{Edwards}, \cite[3.4]{DV}.

\subsection{Theorems of Szpilrajn and V{\"a}is\"{a}l\"{a}}

In  \cite{N}, N\"{o}beling proved that
$$\dim X \leqslant \dim _H X .$$
In \cite[Thm.~2]{Szpilrajn} (see also \cite[Thm. VII~2]{HW}),
Szpilrajn strengthened this showing that
for a separable metric space $X$
$$
\text{ if } m_{n+1}(X) =0, \text{ then }  \dim X \leqslant n.
$$
Moreover, for
$\dim X \leqslant n$
a typical continuous map $f \in C(X,  I^{2n+1})$
satisfies
$m_q (f(X)) = 0$ for any $q >n$,
hence $\dim _H f(X) \leqslant n$
\cite[Thm.~3, Cor.~4]{Szpilrajn}, \cite[Thm.VII~5]{HW}.
The latter fact strengthens the classical 
Lefschetz--Menger--N\"{o}beling--Pontryagin--Tolstowa embedding theorem
\cite[Thm. V~3]{HW}. It also im\-plies that
$\dim X = 
\inf \{ \dim_H Y\},$
where $Y$ runs over all spaces homeomorphic to $X$
\cite[(ii), p.~89]{Szpilrajn}.

For any compactum $X\subset \mathbb R^N$ 
we have (see \cite[Thm. 3.4.3]{DV} and
\cite[6.15]{LV})
$$
\dim X \leqslant \dem X\leqslant \dim_H X .
$$

J.~V{\"a}is\"{a}l\"{a}
proved an analogue of Szpilrajn's theorem for
$\dem $. He showed that
for a compactum
$X\subset \mathbb R^N$
the inequality
$\dem X \leqslant k$
is equivalent to the existence of a homeomorphism
$f:\mathbb R^N\cong\mathbb R^N$
with
$m_{k+1}(f(X)) = 0$.
Moreover,
$\dem X \leqslant \dim _H (f(X))$
 for each homeomorphism
$f:\mathbb R^N\cong\mathbb R^N$,
and this becomes an equality 
for some homeomorphism
$f$ \cite{Vaisala}, \cite[Thm. 3.6.1, 3.6.2]{DV}.
We can also control the size of the isotopy and its 
support \cite[p.~168, Remark]{Vaisala}.

\section{Statements}

All proofs are in Section~\ref{proofs}.

Theorem \ref{main}
is the main result of this paper.
It simultaneously strengthens V{\"a}is\"{a}l\"{a}'s theorem
and provides an answer
to the question posed in Section~\ref{introduction}
about the typical behavior
of projections of compact sets.

On the notion of ``typical'' see Section~\ref{Baire}.

Spaces
$\Homeo_\varepsilon (\overline U,\partial U)$ and
$\Isot _\varepsilon (\overline U,\partial U)$
are defined in Section~\ref{metrics}.

\begin{theorem}\label{main}
Let $X\subset \mathbb R^N$ be a non-empty compactum,
$U$ its bounded open neighborhood, and
 $\varepsilon>0$.
Then

a) a typical element $f$ of the space
$\Homeo _{\varepsilon } (\overline U , \partial U)$
satisfies
$$\dem X = \dim _H f(X)  ,$$

b) a typical element $F = \{ f_t\} $ of the space $\Isot _{\varepsilon } (\overline U , \partial U)$
has the property
$$\dem X = \dim _H f_1( X) .$$
\end{theorem}

This immediately implies

\begin{corollary}\cite{Vaisala}\label{Vai}
For any non-empty compactum $X\subset\mathbb R^N$ and any
$\varepsilon >0$
there exists an
$\varepsilon $-isotopy $F=\{ f_t\} $ of $\mathbb R^N$
with support in
 $O_{\varepsilon } X$
 such that
$\dem X = \dim _H f_1 (X) $.
\end{corollary}

To prove the existence of the desired homeomorphism $f$, V{\"a}is\"{a}l\"{a} uses the equivalence of the inequality $\dem X \leqslant k$ and of the 
ambient embeddability of $X$ into the Menger compactum $M^k_N$
\cite[Thm.~2]{Stanko1971}, see also
\cite[Prop.~1.2]{Edwards}, \cite[Thm.~3.5.1]{DV}. 
V{\"a}is\"{a}l\"{a} achieves 
the control over the isotopy using Edwards's remark
\cite[p. 208--209]{Edwards}: for a given
$\varepsilon $,
instead of the standard Menger com\-pact\-um
one should take a more complicated subset,
depending on $\varepsilon $.
Our work does not use Menger compacta, 
and the proofs rely solely on the definition of the  embedding
dimension $\dem $, standard properties of the Hausdorff measure, and the Baire category.

Now we can find out how often a knot in 
$\mathbb R^3$ has two-dimensional projections,
i.e., how ``common'' the Borsuk-type property is.
Recall that a typical knot in $\mathbb R^3$
is wild \cite{Milnor} and even wild at every
point \cite{Bothe}.

\begin{corollary}\label{case-knots}
Let $\Sigma 
\subset\mathbb R^3$ be a knot or a simple arc,
$U$ its bounded open neighborhood, and
 $\varepsilon>0$.
Then

a) a typical homeomorphism
$f\in \Homeo _\varepsilon (\overline U , \partial U)$ satisfies:
$$
\dim P_\Pi (f (\Sigma )) = \dim_H P_\Pi (f (\Sigma )) =  1
$$
for any non-zero linear subspace
$\Pi \subset \mathbb R^3$;

b) a typical isotopy
$F=\{ f_t \} \in \Isot _\varepsilon (\overline U , \partial U)$ satisfies:
$$
\dim P_\Pi (f_1 (\Sigma )) = \dim_H P_\Pi (f_1 (\Sigma )) =  1
$$
for any non-zero linear subspace
$\Pi \subset \mathbb R^3$.
\end{corollary}

Corollaries \ref{Cantor-tame}, 
\ref{crit2}
include similar statements for Cantor sets.

In the case of zero-dimensional compact sets,
the use of
``double control''
(by the $\varepsilon $ and by the support $U$)
is redundant.

\begin{corollary}\label{Cantor-tame}
Let $K 
\subset\mathbb R^N$ be a Cantor set, 
$N\geqslant 2$. Let
$U$ be a bounded open neighborhood of $K$
such that the intersection of any of its connected components
with $K$ is non-empty.
The following conditions are equivalent:

(i) $K$ is tame;

(ii) there exists a homeomorphism
$f: \mathbb R^N\cong \mathbb R^N
$ with the property:\linebreak
${\dim P_\Pi (f(K) ) =0
}$
for some proper linear  subspace
$\Pi \subset \mathbb R^N$;

(iii)  a typical homeomorphism
$f\in \Homeo  (\overline U , \partial U)$ has the property:
$$
\dim P_\Pi (f (K)) = \dim_H P_\Pi (f (K )) =  0
$$
for any linear subspace
$\Pi \subset \mathbb R^N$;

(iv) a typical isotopy
$F=\{ f_t \} \in \Isot (\overline U , \partial U)$ has the property:
$$
\dim P_\Pi (f_1 (K)) = \dim_H P_\Pi (f_1 (K )) =  0
$$
for any linear subspace
$\Pi \subset \mathbb R^N$;

(v) 
a typical homeomorphism
$f\in \Homeo  (\overline U , \partial U)$
has the property:
the pro\-jec\-ti\-on of $f (K)$
into any non-zero linear subspace is a Cantor set.
\end{corollary}

\begin{corollary}\label{crit2}
Let $K \subset \mathbb R^N$ be a Cantor set, $N\geqslant 3$.
Let
$U$ be a bounded open neighborhood of $K$
such that the intersection of any of its connected components
with $K$ is non-empty.
The following conditions are equivalent:

(i) $K$ is wild;
 
(ii) any homeomorphism
$f: \mathbb R^N\cong \mathbb R^N
$ has the property:
 $\dim P_\Pi (f(K) ) > 0$
for each proper linear subspace
$\Pi \subset \mathbb R^N$;

(iii) any homeomorphism
$f: \mathbb R^N\cong \mathbb R^N
$ has the property:
$\dim P_\Pi (f(K)) = 1$ for each straight line $ \Pi $;

(iv) any homeomorphism
$f: \mathbb R^N\cong \mathbb R^N
$ has the property:
$$\dim P_\Pi (f(K)) \in \{ \dim \Pi , \dim \Pi -1\}$$ for each proper linear subspace
$ \Pi $;

(v) 
a typical homeomorphism
$f\in \Homeo  (\overline U , \partial U)$
has the property:
$$\dim P_\Pi (f(K)) = \dim_H P_\Pi (f(K)) = N-2$$ 
for each $(N-1)$-dimensional linear subspace $\Pi $;

(vi) 
a typical homeomorphism
$f\in \Homeo  (\overline U , \partial U)$
has the property:
 $$\dim _H P_\Pi (f(K)) = \dim \Pi $$ for each $ k\in\{1,\ldots , N-2\} $
 and for almost all (in the sense of Lebesgue measure)
 $k$-dimensional linear subspace%
 \footnote{Here $G(N,k)$ is the Grassmann space
 consisting of all $k$-dimensional linear subspaces of $\mathbb R^N$.}
  $ \Pi \in G(N,k)$;

 (vii)
 a typical isotopy
$F=\{ f_t \} \in \Isot  (\overline U , \partial U)$
has the property:
$$
\dim P_\Pi (f_1(K)) = \dim_H P_\Pi (f_1 (K )) =  N-2
$$
for each $(N-1)$-dimensional linear subspace
$\Pi \subset \mathbb R^N$.
\end{corollary}

\section{Necessary concepts}

\subsection{Baire category}\label{Baire}

\begin{definition}(R.~Baire, 1899)
Let $\mathcal X$ be a non-empty topological space.
A subset $\mathcal A\subset\mathcal  X$ 
is \emph{of first category} in $\mathcal X$
if it can be represented as
$\mathcal A= \bigcup\limits_{i=1}^\infty \mathcal M_i$, where each
$\mathcal M_i$ is nowhere dense in $\mathcal X$;
otherwise,
$\mathcal A$ is \emph{of second category} in $\mathcal X$.
A subset $\mathcal A\subset \mathcal X$ is \emph{resudual} in $\mathcal X$
if
$\mathcal X-\mathcal A$ is of first category in $\mathcal X$.
\end{definition}

For this definition to make sense, we restrict ourselves to considering \emph{Baire} spaces
(i.e., spaces in which
every nonempty open subset has second category in $\mathcal X$;
equivalently, every residual subset is dense in $\mathcal X$.)

A space metrizable by a complete metric
is a Baire space
(R.~Baire for~$\mathbb R$, 1899; F.~Hausdorff, 1914).
A $G_\delta$-subset of a space
metrizable by a complete metric
can also be metrized by a complete metric by P.S.~Alexandroff's theorem
\cite[3.11]{Kechris}.

\begin{definition}
Let $\mathcal X$ be a non-empty Baire space,
$\mathcal A\subset \mathcal X$.
We say that
\emph{a typical element of $\mathcal X$ belongs to $\mathcal A$} 
if the set of elements of $\mathcal X$ that do not belong to $\mathcal A$
has first category.
\end{definition}

\subsection{Spaces of compacta}\label{Baire-hyperspaces}

To illustrate Section~\ref{Baire}, we remind 
some well-known examples;
they will be used below.

\begin{proposition}\label{hyper}
1) Let
$M$ be a complete metric space.
The space 
$\mathcal K (M)$
of all non-empty
compact subsets
of $M$
endowed with the 
Hausdorff metric
is complete and separable \cite[\S 18]{Blaschke}, \cite[4.25]{Kechris}.
The corresponding topology on $\mathcal K ( M)$
coincides with the Vietoris topology
\cite[4.F, 4.21]{Kechris}.

2)  Let
$M$ be a complete metric space without isolated points.
The set $\mathcal C(M)$
of all Cantor sets in $M$
is a dense $G_\delta $ subset of  $\mathcal K (M)$
\cite[Prop.~2]{Kuratowski}.
In particular, $\mathcal C(M)$ is itself a Baire space.
(The Hausdorff metric on
$\mathcal C(M)$
may be non-complete.)

3)
If $M = \mathbb R^N$,
tame Cantor sets in $\mathbb R^N$
form a dense
$G_\delta $-subset $\mathcal T (\mathbb R^N)$
in the space of all Cantor sets
$\mathcal C(\mathbb R^N)$
(\cite[Lemma 1]{Kuzminykh-typ} for the case of $\mathbb R^N$,
or \cite[Lemma 2.1, Remark 2.2]{Gartside}).
\end{proposition}

\subsection{Spaces of homeomorphisms and isotopies}\label{metrics}

Let $(Y, \rho )$be a non-empty metric compactum,
and $(\widetilde Y, \widetilde \rho )$ a metric space.
The space of continuous maps
$C(Y,\widetilde Y)$ is endowed with the metric
$$d(f,g) = \sup\limits _{x\in Y} \rho  (f(x),g(x)).$$
If $\widetilde Y$ is Polish, then
$C(Y,\widetilde Y)$ is also Polish
\cite[(4.19)]{Kechris}. In this paper,
$\widetilde Y$ either coincides with $Y$, or
is a Euclidean space.

We will use the following subspaces of $C(Y,Y)$:

1) $\Homeo (Y)$ --- the set of all homeomorphisms $Y\cong Y$.

2) $\Homeo (Y,A)$ --- the set of all homeomorphisms $f:Y\cong Y$
such that%
\footnote{
Here 
 and below, $\id $ is the identity homeomorphism.
 }
${f|_A = \id }$.

3) $\Homeo _{\varepsilon }(Y,A)$, where $\varepsilon >0$,
is the set of all homeomorphisms $f:Y\cong Y$ such that
$f|_A = \id $ 
and
$d(f,\id ) < \varepsilon $.

For $\varepsilon > \diam Y$
we obviously have
$\Homeo _{\varepsilon }(Y,A)=\Homeo (Y,A)$.

We will use the following subspaces of $C(Y\times I, Y)$:

4) $\Isot (Y)$ is the set of all isotopies of $Y$. 
The distance between
$F=\{ f_t \} ,G=\{ g_t\} \in\Isot (Y)$ is defined by
$$
D(F,G) = 
\sup _{x\in Y , t\in I} \rho (f_t (x) , g_t (x)) .
$$

5) 
$\Isot (Y,A)$ is the set of all isotopies $F:Y\times I\to Y$
such that
$f_t |_A = \id $ for each $t\in I$.

6) 
$\Isot _{\varepsilon } (Y,A)$, where $\varepsilon >0$,
is the set of all isotopies $F:Y\times I\to Y$
such that%
\footnote{Here and below, $\Id $ denotes the ``identical'' isotopy.
That is,
 $F(x,t) = x$ for any $x\in Y$ and $t\in I$.} 
$f_t |_A = \id $ for each $t\in I$
and
$D (F , \Id ) < \varepsilon $.

\begin{proposition}\label{G-delta}
Let $Y$ be a metric compactum,
$A\subset Y$ its closed subset, and
$\varepsilon >0$.

1) The spaces $\Homeo (Y)$, $\Homeo (Y,A)$, 
$\Homeo _{\varepsilon }(Y,A)$
are $G_\delta $-subsets of
$C(Y,Y)$.

2) The spaces $\Isot (Y)$, $\Isot (Y,A)$, 
$\Isot _{\varepsilon }(Y,A)$
are $G_\delta $-subsets of
$C(Y\times I,Y)$.
\end{proposition}

Hence each of the spaces $\Homeo (Y)$, $\Homeo (Y,A)$, 
$\Homeo _{\varepsilon }(Y,A)$, 
$\Isot (Y)$, $\Isot (Y,A)$, 
$\Isot _{\varepsilon }(Y,A)$
is completely mertizable
\cite[4.14]{Kechris}.
When establishing the openness or 
everywhere density of subsets of these spaces, 
we will use 
``standard, uniform'' metrics $d$, $D$ induced from $C(Y,Y)$ and $C(Y\times I,Y)$.
These metrics are in general incomplete.
This is permissible, since the indicated properties
are topological and can be verified in any of the metrics
defining the same topology.

\section{Proofs}\label{proofs}

\subsection{Proof of Proposition~\ref{intro}}\label{intro-proof}

1a) Is a variation of Statement~12
from~\cite{Frolkina-Greece-volume}.
For completeness, we give a proof.

First consider a particular case: $X$ is a Cantor $(N,N-1,N-1)$-set
obtained by the construction of Borsuk
\cite{Borsuk}. 
There is quite a lot of ``freedom'' in this construction,
however, each of these sets is tame
according to the Bing--Keldysh--Osborne tameness criterion
\cite[Thm.~I.4.2]{Keldysh}
(see also
\cite[Stat.~5]{Frolkina-Greece-volume}).

In fact, Borsuk's construction shows that
for any non-empty open subset
$A\subset X$
we have:
$\dim P_\Pi (A) = N-1$ for any
$(N-1)$-plane
$\Pi $. 
Thus in the particular case under consideration we may take
$h_t = \id $ for each $t\in I$.

Now let $X$ be an arbitrary tame Cantor set.
Let 
$\{ L_i\}$ be a defining sequence for $X$ 
\cite[Def.~9]{Frolkina-Greece-volume}. 
There is a positive integer $k$ with the property:
the diameter of each connected component 
$L_k^{(1)},\ldots , L_k^{(s)}$ 
of
$L_k$
is smaller than $\varepsilon $.
We may assume that
the intersection $X_i:= X\cap L_k^{(i)}$
is non-empty for each $i = 1,\ldots , s$
(hence all $X_i$'s are tame Cantor sets).
It is clear that Borsuk's construction
provides Cantor sets of arbitrarily small diameter
(apply a homothety is necessary).
Hence 
for each $i = 1,\ldots , s$
there exists
a tame Cantor $(N,N-1, N-1)$-set
$Y_i \subset \Int L_k^{(i)}$ 
\cite{Borsuk}.
For any non-empty open subset
 $A\subset Y_1\cup
\ldots \cup Y_s$
the projection of 
$A$ into any $(N-1)$-plane
has dimension
$(N-1)$.
The proof of Theorem~I.4.2 from
\cite{Keldysh}
(see also \cite[Stat.~6]{Frolkina-Greece-volume})
implies:
there is an isotopy $\{ h_t \} : \mathbb R^N\cong \mathbb R^N$
such that
$h_t |_{\mathbb R^N - \cup _{i=1}^s L_k^{(i)} } = \id $ 
for each $t\in I$,
and $h_1(X_i) = Y_i$ for each $i=1,\ldots , s$.
This is the desired isotopy.

1b) is proved similarly to 1a).
Now $Y_i$'s are
choosed differently:
they should be tame
Cantor sets in general position
with respect to all projections
\cite[Thm.~5]{Cobb-projections},
\cite{Frolkina-arxiv}, \cite{Bogatyi}.
The union
$Y:= Y_1\cup \ldots \cup Y_s$
is a Cantor set \cite[(7.4)]{Kechris}
and it is tame
\cite[Thm.~6.1]{Bing-tame},
\cite[Thm.~8]{Osborne1}.
Moreover,  $Y_i$'s can be taken so that
 $Y$ has general position with respect to all projections.

2) 
is proved similarly to Theorem~2 from
\cite{Frolkina-Greece-volume}.

3) Any non-countable compactum $X$ contains a Cantor set
$K$ \cite[(6.2)]{Kechris}.
Moreover: if $X$ is embedded in $\mathbb R^N$,
then
 $K$ can be chosen to be tame in $\mathbb R^N$.
 (For this, use the tameness criterion
\cite[Thm.~I.4.2]{Keldysh},
\cite[Stat.~5]{Frolkina-Greece-volume}.)
To finish the proof, apply~1a).

4) By \cite[Satz~4]{Bothe-U},
$\dem X =1$.
Together with
\cite{Vaisala}
this implies:
there is an
$\varepsilon $-isotopy $\{ h_t \} : \mathbb R^3 \cong \mathbb R^3$
with support in $O_\varepsilon X$ such that
$\dim _H h_1( X) = 1$.
Orthogonal projection can not raise Hausdorff dimension;
hence
$\dim _H  P_\Pi (h_1( X) ) \leqslant  1$ for any linear subspace
$\Pi \subset \mathbb R^3$.
For a non-zero subspace
$\Pi \subset \mathbb R^3 $
the equality $\dim _H P_\Pi (h_1 (X)) = 0$
is possible only in two cases:

(i) $\dim \Pi = 2$, and $h_1(X)$ 
is a segment perpendicular to $\Pi$;

(ii) $\dim \Pi = 1$, and $h_1(X)$ is contained in a plane 
perpendicular to $\Pi $.

In both cases, it is easy to
``tweak'' the isotopy $\{ h_t \}$
to obtain the desired property.

\subsection{Proof of Proposition~\ref{G-delta}}

1) It is known that
$\Homeo (Y)$ 
is a  $G_\delta $-subset of
$C(Y,Y)$  \cite[Lemma 2.2.1]{DV}.
Observe that
$$
\Homeo (Y,A) =
\Homeo (Y) \cap
\{ f\in C(Y,Y) \ : \ f|_A = \id 
\}
.
$$
The set
$\{ f\in C(Y,Y) \ : \ f|_A = \id 
\}$ 
is closed in $C(Y,Y)$,
hence is a $G_\delta $-subset of $C(Y,Y)$.
Consequently
$\Homeo (Y,A)$ 
is a $G_\delta $-subset of
$C(Y,Y)$.

We further observe that
$$
\Homeo _{\varepsilon } (Y,A)
=
 \Homeo (Y,A) \cap
\{ f\in  C(Y,Y) \ : \ d (f, \id ) < \varepsilon \}  .
$$
The second of the two intersected sets is open in $C(Y,Y)$.
This proves 1).
Part 2) can be proved by similar reasoning.

\subsection{Proof of Theorem~\ref{main}}

The proof is very close to the arguments of
S.~Eilenberg,
see the proof of Theorem~3 in \cite{Szpilrajn}
or the proof of Theorem VII~5 in \cite{HW}.

a) For $q>0$
and $k\in\mathbb N$
define
$W_{q ,  k }$ as the set of all
$f\in \Homeo _\varepsilon  (\overline U , \partial U)$
with the property:
for some finite collection
$V_1, \ldots ,V_s$
of open subsets of $\mathbb R^N$
we have
$$
f(X) \subset V_1 \cup \ldots \cup V_s \subset U 
\quad
\text{ and }
\quad 
\sum_{i=1}^{s} (\diam V_i)^q < \frac{1}{k}  .
$$
Evidently each
$W_{q ,  k }$ is open in $\Homeo _\varepsilon (\overline U , \partial U)$.
Hence
the set
$$T := 
 \bigcap\limits_{\substack{q= \dem X +\frac{1}{n}\\ n\in \mathbb N}}
\left( \bigcap\limits_{k \in \mathbb N} W_{q,  k } \right) 
$$
is a $G_\delta $-subset of $\Homeo _\varepsilon (\overline U , \partial U)$.

\emph{For each
$f\in T$
we have $\dim _H  f(X) = \dem X $.
}
In fact, take any
$f\in T$.
For each $q >\dem X $ we have ${m_q (f(X)) = 0}$ 
by
\cite[VII (1D)]{HW} and by definition of $T$.
Hence
$\dim _H  f(X) \leqslant \dem X$.
Recall  that $\dem $ is invarint under homeomorphisms of $\mathbb R^N$.
We get
$\dem X = \dem f(X)$ \cite[Thm.~3]{Stanko1970}.
This together with the inequality
$\dem f(X) \leqslant \dim _H  f(X)$
\cite[6.15]{LV}, \cite[Thm.~3.6.2]{DV}
implies
$\dim _H  f(X) = \dem X $.

\emph{Fix any $q >\dem X $ and any
$k \in\mathbb N$. Let us show that the set $W_{q ,  k }$ 
is dense in $\Homeo _\varepsilon (\overline U , \partial U)$.
}
To show this, take any
$g\in \Homeo _\varepsilon  (\overline U , \partial U)$
and $\gamma >0$.
Let us construct 
an $f\in W_{q,   k }$ with $d(f,g) < \gamma $.
Take a positive number
$$
\delta < 
\min \Bigl\{ 
\gamma ;
\quad 
\frac13 d (g(X) , \partial U); 
\quad
\varepsilon -  d(g,\id )
\Bigr\} .
$$
Since $\dem g(X) = \dem X $, there exists
a finite polyhedron $P\subset \mathbb R^N$
of dimension $\dim P = \dem X $
and a
$\delta $-pseudoisotopy $\Phi = \{ \varphi _t \} : \mathbb R^N\to \mathbb R^N$
with support in $O_{\delta } (g(X)) $
such that
$\varphi _1 (g(X)) = P$ \cite[Def.~1]{Stanko1970}.
Note that $O_{\delta } P \subset O_{2\delta } (g(X)) \subset U $.

We have $m_q P = 0$ since $q >\dem X = \dim P $.
Hence there is a representation
$P=L_1 \cup \ldots \cup L_r$ with
$$
\sum\limits_{i=1}^{r} (\diam L_i)^q  < \frac{1}{k}
.
$$
Let $\widetilde L_1, \ldots , \widetilde L_r$ be open
subsets of
$\mathbb R^N$
with
$$
\sum \limits_{i=1}^r (\diam \widetilde L_i)^q  < \frac{1}{k},
\quad
\bigcup \limits_{i=1}^r \widetilde L_i \subset U,
\quad
L_i \subset \widetilde L_i
\text{ for each } 
i=1,\ldots ,r.$$
Continuity of $\Phi : \mathbb R^N \times I\to\mathbb R^N$
and compactness of $g(X)$ imply:
$\varphi _\tau (g(X)) \subset \bigcup\limits _{i=1}^r \widetilde L_i$
for some $\tau \in [0,1)$.

Define $f:= \varphi _\tau \circ g$.
Then
$$
d(f,g) = d(\varphi _\tau \circ g ,g) < d(\varphi _\tau , \id ) \leqslant \delta <\gamma 
$$
and
$$
d(f,\id )\leqslant d(f,g)+d(g,\id ) <  \delta
+ \left( \varepsilon  - \delta \right) 
=\varepsilon .
$$
Together with
$$
f(X) = \varphi _\tau (g(X)) \subset \bigcup\limits _{i=1}^r \widetilde L_i \subset U 
\quad
\text{ and }
\quad
\sum \limits_{i=1}^r (\diam \widetilde L_i)^q < \frac{1}{k} 
$$
this implies: $f\in W_{q ,  k }$. Thus a) is proved.

Part b) is proved by analogous reasoning.
For $q >0$ and $k\in\mathbb N$
instead of $W_{q,  k}$
consider
$\mathbb W_{q ,  k  }$ consisting of all
$F\in \Isot _\varepsilon (\overline U, \partial U)$
with the property:
for some finite collection
$V_1, \ldots ,V_s$
of open subsets of $\mathbb R^N$
we have
$$
f_1 (X) \subset V_1 \cup \ldots \cup V_s \subset U
\quad
\text{ and }
\quad\sum_{i=1}^{s} (\diam V_i)^q  < \frac{1}{k} 
.
$$
We similarly define
$$\mathbb T := 
 \bigcap\limits_{\substack{q = \dem X+\frac{1}{n}\\ n\in \mathbb N}}
\left( \bigcap\limits_{k \in \mathbb N} \mathbb W_{q ,  k } \right) ,
$$
this is a $G_\delta $-subset of
 $\Isot _\varepsilon (\overline U, \partial U)$.

The density of $\mathbb W_{q, k}$ in
$\Isot _\varepsilon (\overline U, \partial U)$ for $q >\dem X$ and $k \in \mathbb N$
is somewhat more involved than for part a),
but is not difficult.
Briefly, given an isotopy
$G=\{ g_t \} \in \Isot _\varepsilon (\overline U, \partial U)$ 
and a number $\gamma >0$,
then constructing
an isotopy $F=\{ f_t \} \in \mathbb W_{q, k}$ satisfying
$D(F,G) < \gamma $ can be done as follows:
successively apply the isotopy $G$,
and then some isotopy obtained
by incompletely applying a small pseudoisotopy,
mapping $g_1 (X) $ to a polyhedron of dimension equal to $\dem X$.

\subsection{Proof of Corollary~\ref{case-knots}}

a) 
Consider two cases.

{\it A particular case.} The set $\Sigma \subset \mathbb R^3$
is embedded equivalently to a standard circle
in $\mathbb R^2 \times \{ 0\}$ or a segment
in $\mathbb R \times \{ 0\}\times \{ 0\}$.
In this case, a) is easy.
Moreover:
the set of all homeomorphisms
$f\in \Homeo _\varepsilon (\overline U , \partial U)$ such that
$f (\Sigma )$ lies in a $2$-plane
is nowhere dense in
$\Homeo _\varepsilon (\overline U , \partial U)$.

{\it The general case.}  
The embedding $\Sigma \subset \mathbb R^3$
is non-equivalent to a standard circle
or a segment.
Obviously $\dim P_{\Pi } \Sigma \neq 0 $
for any $2$-plane $\Pi $.

Let us show that 
$\dim P_{\ell } \Sigma \neq 0$
for any straight line $\ell $.
Suppose the contrary.
Then $\Sigma  \subset L$, where $L$ is a $2$-plane
perpendicular to $\ell $.
By Jordan--Sch{\"o}nflies Theorem \cite[Thm.~10.4]{Moise},
there is a homeomorphism $h:L\cong L$
such that $h(\Sigma )$ is a circle or a segment.
Hence $H= h\times \id : L\times \ell \cong L\times \ell $
is a homeomorphism of $\mathbb R^3$
which takes $\Sigma $ into a circle or a segment, a contradiction.

Recall that
$\dem \Sigma = 1$ 
for any knot and any simple arc $\Sigma \subset \mathbb R^3$
\cite[Satz~4]{Bothe-U}.
This together with Theorem~\ref{main} implies:
$\dim _H  f(\Sigma ) = 1$
for a typical homeo\-mor\-phism
$f\in \Homeo _\varepsilon (\overline U , \partial U)$.
For any $2$-plane or line
$\Pi \subset \mathbb R^3$
we then get
$$
0\neq \dim P_\Pi (f(\Sigma )) \leqslant 
\dim _H  P_\Pi (f(\Sigma )) \leqslant \dim _H  f(\Sigma ) = 1.
$$
Hence $\dim P_\Pi (f(\Sigma )) =1$, and a) is proved.

Similar arguments prove Part~b).

\subsection{Proof of Corollary~\ref{Cantor-tame}}

$(i) \Rightarrow (iii)$. 
For any tame Cantor set $K\subset \mathbb R^N$
we have
$\dem K = 0$.
By Theorem~\ref{main},
$\dim _H   f(X)    =0 $
for a typical element $f\in
\Homeo  (\overline U , \partial U)$.
Hausdorff dimension can not raise under projections.
Hence
$$
\dim  P_\Pi (f (K))  \leqslant \dim _H   P_\Pi (f (K)) 
\leqslant \dim _H  f (K)  = 0
$$
for any linear subspace 
$\Pi \subset \mathbb R^N$.

$(i) \Rightarrow (iv)$ is proved by similar arguments.

It is evident that any of 
$(iii), (iv), (v)$ implies $(ii)$.
 
Let us prove that $(ii)\Rightarrow (i)$.

\emph{Case 1.} For $N=2$ there is nothing to prove.

\emph{Case 2.}
Let $N=3$,
and 
$\dim P_\Pi (f(K)) = 0$
for some $2$-plane $\Pi $.
We may assume that $\Pi = \mathbb R 
\times \mathbb R \times \{ 0 \}$ (a coordinate plane).
By Antoine's theorem,
any zero-dimensional compactum
in plane is tame
\cite[{\bf 75}, p.~87--89]{Antoine-diss},
\cite[Cor. II.3.2, II.3.3]{Keldysh},
\cite[Chap.~13]{Moise}.
Hence there is a homeomorphism $h:\Pi \cong \Pi $
with
$h( P_\Pi (f(K))) \subset \mathbb R \times \{ 0\}$.
Define a homeomorphism
$\widetilde h = h\times \id _\mathbb R : 
\mathbb R^3 \cong\mathbb R^3$,
we have
$\widetilde h (f(K)) \subset \mathbb R \times \{0\}\times \mathbb R$.
Applying Antoine's theorem again, we get: $f(K)$ is tame in $\mathbb R^3$.
Consequently $K$ is tame.

\emph{Case 3.}
$N\geqslant 3$, and 
$\dim P_\Pi (f(K)) = 0$
for some straight line $\Pi $.
Necessary statement follows from
\cite[Cor.~6.2]{Frolkina-Greece-volume}.

\emph{Case 4 (the general one).}
$N\geqslant 4$, 
and $\dim P_\Pi (f(K)) = 0$
for some linear subspace $\Pi $ with
 $2\leqslant \dim \Pi \leqslant N-1$.
Let $\widetilde \Pi  \subset \mathbb R^N$ be an $(N-1)$-dimensional linear 
subspace
containing $\Pi $.
We have
$\dim P_{\widetilde \Pi } (f(K)) \leqslant N-1-\dim \Pi  \leqslant N-3$,
and desired statement follows from
\cite[Thm. 5.3, 5.4]{Wright}, \cite{Walsh-Wright}.
(These two papers use the fact:
a Cantor set $X\subset \mathbb R^N$ is tame
iff $\mathbb R^N - X$ is
$1$-LC, see \cite[Thm. 3.4.11, Exer. 3.4.1]{DV}. 
For the cases $\dim \Pi = N-1 $ or $\dim \Pi = N-2$,
a more direct proof of tameness of $K$
is given in
\cite[Thm.~3]{Frolkina-Greece-volume},
where Bing--Keldysh--Osborne tameness criterion is used.)

Finally, we show that
$(i)\Rightarrow (v)$.
For a non-empty open set
$U\subset \mathbb R^N$
by $\mathcal T'(U)$
denote the space of all tame Cantor 
sets
$Y\subset \mathbb R^N$
which have non-empty intersection with any
connected component of $U$. The space
$\mathcal T'(U)$ is endowed with the Hausdorff metric.
It is completely metrizable,
see Proposition~\ref{hyper}.
By
\cite[(4.19)]{Kechris} together with Proposition~\ref{G-delta},
the space $\Homeo ( \overline U , \partial U)$
is also completely mertizable.
The proof in
\cite[I.4.2]{Keldysh}
(see also
\cite[Stat.~6]{Frolkina-Greece-volume})
implies that the group action
$$
\Psi : 
\Homeo ( \overline U , \partial U)
\times 
\mathcal T'(U)
\to
\mathcal T'(U),
\quad
(h,X)\mapsto h(X) .
$$
is transitive.
By Effros  Theorem \cite{Effros}, \cite{Ancel-Effros}, \cite{vanMill-Effros}, the map
$$
\Psi _K : 
\Homeo ( \overline U , \partial U)
\to
\mathcal T'(U),
\quad
h\mapsto h(K) 
$$
is an open surjection.
Now $(v)$
follows since a typical element
of $\mathcal T'(U)$ is a Cantor set having
general position with respect to all projections
\cite{Frolkina-arxiv}, \cite{Bogatyi}.

\subsection{Proof of Corollary \ref{crit2}}

Implications
$(v), (vi), (vii)\Rightarrow (i)\iff (ii)$ follow from Corollary~\ref{Cantor-tame}.
Implications
$(ii)\Rightarrow (iii)\Rightarrow (i)$
and
$(iv)\Rightarrow (i)$
are evident.
(For $(i)\Rightarrow (iii)$ see also
\cite[Thm.~3]{Frolkina-Greece-volume}.)

$(i)\Rightarrow (iv)$.
For
 $N=3$ see the proof of Corollary~\ref{Cantor-tame}, namely,
of $(ii)\Rightarrow (i)$.
Now let $N\geqslant 4$.
Suppose the contrary: for some homeomorphism
$f:\mathbb R^N\cong\mathbb R^N$
and some proper
linear subspace $\Pi $ we have
$\dim P_\Pi (f(K)) \leqslant \dim \Pi -2$.
Take an $(N-1)$-dimensional linear subspace 
$\tilde \Pi \supset \Pi $.
Then
$$\dim P_\Pi (f(K)) \leqslant ( \dim \Pi -2) + (\dim \tilde \Pi - \dim \Pi )
= N-3.$$
By \cite{Wright}, \cite{Walsh-Wright}, $f(K)$ is tame.
Consequently $K$ is tame.

Let us prove
$(i)\Rightarrow (v)$.
For any wild Cantor set $K\subset \mathbb R^N$
we have
$\dem K = N-2 $.
Hence 
$\dim _H  f(X)   =N-2 $
for a typical element $f\in
\Homeo  (\overline U , \partial U)$.
Then for any $(N-1)$-plane
$\Pi \subset \mathbb R^N$ we have
$$
\dim P_\Pi (f (K)) \leqslant \dim _H P_\Pi (f (K))
\leqslant \dim _H (f (K)) = N-2 .
$$

Suppose that
$\dim P_\Pi (f (K)) \leqslant N-3$
 for some $(N-1)$-plane
$\Pi $.
Then $f(K)$ is tame
\cite[Thm. 5.4]{Wright}, \cite[Thm. 3]{Frolkina-Greece-volume}.
Consequently $K$ is tame, a contradiction.
Thus $\dim P_\Pi (f (K)) = N-2$
for any $(N-1)$-plane
$\Pi $.

The proof of
$(i)\Rightarrow (vii)$
is similar to that of $(i)\Rightarrow (v)$.

Finally, we show that
$(i)\Rightarrow (vi)$.
As above:
for a typical element $f\in
\Homeo  (\overline U , \partial U)$
we have
$\dim _H  f(X)   =N-2 $.
To get necessary statement, apply
\cite[Thm.~6.2]{Falconer}
(see p.~105--106 of the cited book for historic details and references).

\end{document}